%
%
%
%
\font\tenbf=cmbx10

\font\eightrm=cmr8
\font\eightit=cmti8
\font\germ=eufm10

\def\s{\hbox{\germ S}}

\def\sectiontitle#1\par{\vskip0pt plus.1\vsize\penalty-250
\vskip0pt plus-.1\vsize\bigskip\vskip\parskip
\message{#1}\leftline{\tenbf#1}\nobreak\vglue 5pt}
\def\wt{\widetilde}
\def\eno{\eqalignno}

\def\mod{\hbox{\rm mod~}}

\def\ds{\displaystyle}
\magnification=\magstep1
\parindent=15pt
\nopagenumbers
\baselineskip=10pt
\baselineskip=13pt
\headline{\ifnum\pageno=1\hfil\else%
{\ifodd\pageno\rightheadline \else \leftheadline\fi}\fi}
\def\rightheadline{\hfil\eightit 
Quantum Grothendieck polynomials
\quad\eightrm\folio}
\def\leftheadline{\eightrm\folio\quad 
\eightit 
Anatol N. Kirillov   
\hfil}
\voffset=2\baselineskip
\centerline{\tenbf 
QUANTUM\hskip 0.1cm GROTHENDIECK\hskip 0.1cm POLYNOMIALS 
}
\vglue 24pt
\centerline{\eightrm 
ANATOL N. KIRILLOV
}
\baselineskip=12pt
\centerline{\eightit 
Centre de Recherches Mathematiques, Universite de Montreal,
}
\baselineskip=10pt
\centerline{\eightit 
C.P. 6128, Succursale, Centre-Ville, Montreal (Quebec), H3C 3J7, Canada
}
\baselineskip=10pt
\centerline{\eightit
and
}
\baselineskip=10pt
\centerline{\eightit  
Steklov Mathematical Institute,
}
\baselineskip=10pt
\centerline{\eightit 
Fontanka 27, St.Petersburg, 191011, Russia
}
\vglue 15pt
\centerline{\eightrm ABSTRACT}
{\rightskip=1.5pc
\leftskip=1.5pc
\eightrm\parindent=1pc 
We study the algebraic aspects of (small) quantum equivariant $K$-theory
of flag manifold. Lascoux-Sch\"utzenberger's type formula for
quantum double and quantum double dual Grothendieck polynomials and the 
quantum Cauchy identity for quantum Grothendieck polynomials 
are obtained. 
} 
\vglue12pt
\baselineskip=13pt
\overfullrule=0pt
\def\qed{\hfill$\vrule height 2.5mm width 2.5mm depth 0mm$}
\vskip 0.5cm

{\bf \S 1. Introduction.}
\vskip 0.4cm

There exists a remarkable hierarchy of polynomials related with the
generalized cohomology theories of the flag manifold:

$\bullet$ double and double dual Grothendieck polynomials 
$G_w^{(\beta )}(x,y)$ and ${\cal H}_w^{(\beta )}(x,y)$, corresponding 
to the equivariant $K$-theory;

$\bullet$ double Schubert polynomials $\s_w(x,y)=
G_w^{(\beta )}(x,y)|_{\beta =0}$, corresponding to the equivariant
cohomology theory;

$\bullet$ Schubert polynomials $\s_w(x)=\s_w(x,y)|_{y=0}$, corresponding
to the singular cohomology theory.

The Grothendieck polynomials, which are the $K$-theoretic analog of
Schubert polynomials, was introduced by Alexander Grothendieck in his
study of the general Riemann-Roch Theorem [G]. Algebraic definition
and fundamental properties of Grothendieck and Schubert polynomials were
introduced and developed in series of Alain Lascoux and 
Marcel-Paul Sch\"utzenberger papers, see, for example, [LS], [L],
for more references see [M].

The hierarchy, mentioned above, reflects the well-known fact, 
that the corresponding generalized cohomology theories 
are the only ones which have a polynomial group low.

In this paper we continue the quantization of some interesting classes of
polynomial 
and introduce the quantum double, $\wt G_w(x,y)$,
and quantum double dual, $\wt{\cal H}_w(x,y)$, Grothendieck polynomials and
investigate their properties.

Quantum double Schubert polynomials, $\wt{\s}_w(x,y)$, were introduced
and their properties were studied in [KM].

Quantum Schubert polynomials, $\s_w^q(x)$, were introduced and their
properties were studied in [FGP]. Another approach to quantum Schubert 
polynomials 
is considered in [KM].

The structure of this paper is the following. In \S 2 some results on double 
Schubert and double Grothendieck polynomials,
quantum analogs of which we would like to construct, are collected. Probably,
the formulae (5),(6),(7) was not appear in earlier publications.
In \S 3 we recall some results from [FGP], which will be used in \S 4.
In \S 4 we introduce quantum double and quantum double dual Grothendieck
polynomials and prove the quantum Cauchy identity.

{\it Acknowledgement.} I would like to thank Dr.~T.~Maeno for fruitful
discussions and collaboration.

\vskip 0.4cm

{\bf \S 2. Classical Grothendieck polynomials.}
\bigbreak

In this section we give a brief review of the theory of Grothendieck and
Schubert polynomials created by A.~Lascoux and M.P.-Sch\"utzenberger,
see e.g. [L], [LS] and [M]. A proof of Cauchy's formula for Grothendieck 
polynomials can be found in [FK], and that of Pieri's rule for double 
Grothendieck polynomials can be obtained using the generalization of
method from [KV]. 
\bigbreak

{\bf 2.1. Divided and isobaric divided differencies.}
\bigbreak

Let $x_1,\ldots x_n$ and $\beta$ be independent set of variables, and let 
$P_n:=P_{n,\beta}={\bf Z}[x_1,\ldots ,x_n,\beta ]$ for each $n\ge 1$.

Let us denote by $\Lambda_n:=\Lambda_{n,\beta }={\bf Z}{[\beta ]}[x_1,
\ldots , x_n]^{S_n}\subset P_n$ the subring of symmetric polynomials in 
$x_1,\ldots ,x_n$ with polynomial coefficients and by
$$H_n:=H_{n,\beta}=\left\{\sum_{I=(i_1,\ldots ,i_n)}a_Ix^I~|~a_I\in{\bf Z}
{[\beta ]},~0\le i_k\le n-k,~\forall k\right\}
$$
the additive subgroup of $P_n$ spanned by all monomials $a_Ix^I$ with 
$a_I\in{\bf Z}{[\beta ]}$ and
$I\subset\delta :=\delta_n=(n-1,n-2,\ldots ,1,0)$. 

For $i\le i\le n-1$ let us define a $\Lambda_n$-linear operator $\partial_i$
acting on $P_n$
$$(\partial_if)(x)={f(x_1,\ldots ,x_i,x_{i+1},\ldots ,x_n)-f(x_1,\ldots ,
x_{i+1},x_i,\ldots ,x_n)\over x_i-x_{i+1}}
$$
Divided difference operators $\partial_i$ satisfy the following relations
$$\eno{
&\partial_i^2=0,\cr 
&\partial_i\partial_j=\partial_j\partial_i , \ {\rm if}
\ |i-j|\ge 2, & (1)\cr 
&\partial_i\partial_{i+1}\partial_i=\partial_{i+1}
\partial_i\partial_{i+1},}
$$
and the Leibnitz rule
$$\partial_i(fg)=\partial_i(f)g+s_i(f)\partial_i(g).
$$
Let us also introduce the operators $\pi_i:=\pi^+_i$ and $\pi_i^-$ 
($1\le i\le n-1$) defined by
$$(\pi_i^{\pm}f)(x)=(\partial_if)(x)\pm \beta\partial_i(x_{i+1}f(x)).
$$
In place of (1) we have
$$\eno{
& (\pi_i^{\pm})^2=\pm\beta\pi_i,\cr
&\pi_i\pi_j=\pi_j\pi_i, \ {\rm if} \ |i-j|\ge 2,\cr
&\pi_i\pi_{i+1}\pi_i=\pi_{i+1}\pi_i\pi_{i+1},}
$$
and the modified Leibnitz rule
$$\pi_i(fg)=\pi_i(f)g+s_i(f)(\pi_i+\beta )(g),
$$
It is clear that if $f\in\Lambda_n$, then $(\pi_i+\beta )(f)=0$, and
$\pi_i(fg)=f\pi_i(g)$.

For any permutation $w\in S_n$. let us denote by $R(w)$ the set of reduced
decompositions for $w$, i.e. sequences $(a_1,\ldots ,a_p)$ such that
$w=s_{a_1}\ldots s_{a_p}$, where $p=l(w)$ is the length of the permutation 
$w\in S_n$, and $s_i=(i,i+1)$ is the simple transposition that interchanges 
$i$ and $i+1$.

For any sequence ${\bf a}=(a_1,\ldots ,a_p)$ of positive integers, we define
$$\partial_{\bf a}=\partial_{a_1}\ldots\partial_{a_p},\ \pi_{\bf a}=
\pi_{a_1}\ldots\pi_{a_p}.
$$

{\bf Proposition 1} ([M], (2.5), (2.6), (2.15)).
{\it If ${\bf a}, {\bf b}\in R(w)$, then $\partial_{\bf a}=\partial_{\bf b}$
and $\pi_{\bf a}=\pi_{\bf b}$.}

If ${\bf a}$ is not reduced, then $\partial_{\bf a}=0$. 
It follows from Proposition 1, that the operators 
$\partial_w=\partial_{\bf a}$ and $\pi_w
=\pi_{\bf a}$ are well-defined, where ${\bf a}$ is any reduced word for $w$.

\bigbreak
{\bf 2.2. Double Grothendieck and double dual Grothendieck polynomials.}
\medbreak

Let $x=(x_1,\ldots ,x_n)$, $y=(y_1,\ldots ,y_n)$ be two sets of independent
variables and
$$G_{w_0}(x,y):=\prod_{i+j\le n}(x_i+y_i).
$$

{\bf Definition 1} (Lascoux-Sch\"utzenberger [LS],[L], see also [FK]). 
{\it For each 
permutation $w\in S_n$ the double Grothendieck polynomial $G_w(x,y)$ and
double dual Grothendieck polynomial ${\cal H}_w(x,y)$ are defined to be
$$\eno{
G_w(x,y)&=\pi_{w^{-1}w_0}G_{w_0}(x,y),\cr
{\cal H}_w(x,y)&=\psi_{w^{-1}w_0}G_{w_0}(x,y),}
$$
where $w_0$ is the longest element of $S_n$, and $\psi_w:=\ds\sum_{v\le w}
\beta^{l(w)-l(v)}\pi_v$.}

\bigbreak
{\bf Proposition 2} ([L], [LS]).

$\bullet$ {\it (Stability) Let $m>n$ and let $i: ~S_n\hookrightarrow S_m$
be the embedding. Then
$$G_w(x,y)=G_{i(w)}(x,y).
$$

$\bullet$ $\pi_w=\sum_{v\le w}(-\beta )^{l(w)-l(v)}\psi_v$, and hence,
$$\eno{
&{\cal H}_w(x,y)=\sum_{w\le v}\beta^{l(v)-l(w)}G_v(x,y)\cr
&G_w(x,y)=\sum_{w\le v}(-\beta )^{l(v)-l(w)}{\cal H}_v(x,y).}
$$

$\bullet$ Grothendieck and dual Grothendieck polynomials $G_w(x):=G_w(x,0)$
and ${\cal H}_w(x):={\cal H}_w(x,0)$ from the ${\bf Z}[\beta ]$-bases of
$H_n$.

$\bullet$ (Cauchy formula, [FK])
$$\sum_{w\in S_n}{\cal H}_w(x,\ominus z)G_{ww_0}(y,z)=\prod_{i+j\le n}(x_i+y_i+
\beta x_iy_i),\eqno (2)
$$
where $\ominus z:=\ds\left({-z_1\over 1-\beta z_1},{-z_2\over 1-\beta z_2},
\ldots ,{-z_n\over 1-\beta z_n}\right)$.

\medbreak
$\bullet$ ${\cal H}_w^{(\beta )}(x,y)=G_{w^{-1}}^{(-\beta)}(y,x)$, $w\in S_n$.}

Hence, $(\pi_{ww_0})^{(y)}G_{w_0}(x,y)=G_{w^{-1}}^{(-\beta )}(y,x)=
{\cal H}_w^{(\beta )}(x,y)$.

\bigbreak
{\bf Examples.} 1) Double Grothendieck and double dual Grothendieck 
polynomials for $S_3$
$$\eno{
&G_{121}(x,y)=(x_1+y_1)(x_1+y_2)(x_2+y_1),\cr
&G_{12}(x,y)=(x_1+y_1)(x_2+y_1)(1-\beta y_2),\cr
&G_{21}(x,y)=(x_1+y_1)(x_1+y_2)(1-\beta y_1),\cr
&G_1(x,y)=(x_1+y_1)(1-\beta y_1)(1-\beta y_2),\cr
&G_2(x,y)=(x_1+x_2+y_1+y_2+\beta x_1x_2-\beta y_1y_2)(1-\beta y_1),\cr
&G_{{\rm id}}(x,y)=(1-\beta y_1)^2(1-\beta y_2);\cr
&{\cal H}_{121}(x,y)=G_{121}(x,y),\cr
&{\cal H}_{12}(x,y)=(x_1+y_1)(x_2+y_1)(1+\beta x_1),\cr
&{\cal H}_{21}(x,y)=(x_1+y_1)(x_1+y_2)(1+\beta x_2),\cr
&{\cal H}_1(x,y)=(x_1+y_1)(1+\beta x_1)(1+\beta x_2),\cr
&{\cal H}_2(x,y)=(x_1+x_2+y_1+y_2+\beta x_1x_2-\beta y_1y_2)(1+\beta x_1),\cr
&{\cal H}_{\rm id}(x,y)=(1+\beta x_1)^2(1+\beta x_2).}
$$
2) Let $w\in S_n$ be a dominant permutation with code $c(w)=(c_1\ge c_2
\ge\cdots\ge c_n\ge 0)$. Then
$${G_w(x,y)\over G_{\rm id}(x,y)}=\prod_{k=1}^n\prod_{i=1}^{c_k}{x_k+y_i\over
1-\beta y_i}.
$$
\bigbreak

{\bf 2.3. Scalar product.}
\medbreak

Let us define a scalar product on $P_{n.\beta}$ with values in $\Lambda_{n\beta}$,
by the rule
$$\langle\langle f,g\rangle\rangle =\pi_{w_0}(f,g),~~
f,g\in P_{n,\beta},
$$
where $w_0$ is the longest element in $S_n$. The scalar product 
$\langle\langle ,\rangle\rangle$ defines a non-degenerate 
pairing $\langle\langle ,\rangle\rangle_0$ on the quotient ring
$P_{n,\beta}/I_{n,\beta}$, where $I_n:=I_{n,\beta}$ is the ideal in 
$P_{n,\beta}$ generated by the elementary symmetric polynomials $e_1(x),
\ldots ,e_n(x)$.

\bigbreak
{\bf Proposition 3.} {\it

$\bullet$ $\langle\langle \pi_wf,g\rangle\rangle =
\langle\langle f,\pi_{w^{-1}}g\rangle\rangle ,~~
f,g\in P_{n,\beta}$.

$\bullet$ $\langle\langle f,G_w\rangle\rangle_0=
\eta (\pi_{w_0w}f)$, where $\eta :P_{n,\beta }\to{\bf Z}[\beta ]$ is the 
homomorphism defined by $\eta (1)=1$, $\eta (x_i)=0$, $1\le i\le n$.

$\bullet$ (Interpolation formula) If $f\in H_n$, then
$$f(x)G_{\rm id}(x,-y)=\sum_{w\in S_n}{\cal H}_w(x,-y)\pi_w^{(y)}f(y).
\eqno (3)
$$

$\bullet$ (Orthogonality, [LS], [L]) Let $u,v\in S_n$, then}
$$\langle\langle {\cal H}_u(x),G_v(x)\rangle\rangle_0=\cases{1,& if 
$v=w_0u$,\cr 0,& otherwise}
$$

As for every $w\in S_n$, $\langle\langle G_w(x,y),1\rangle
\rangle =\pi_{w_0}(G_w(x,y))=(-\beta )^{l(w_0w)}G_{\rm id}(x,y)$,
the fact that
$$\langle\langle {\cal H}_w(x,y),1\rangle\rangle =\cases{G_{\rm id}(x,y),& if 
$w=w_0$,\cr 0,& otherwise}
$$
generalizes the property of the Moebius function (for the Ehresman-Bruhat order
on $S_n$) to be equal to $\pm 1$, cf. [LS].

\bigbreak
{\bf Remark.} We have
$$\eno{
G_{\rm id}(x,y)&=\prod_{k=1}^{n-1}(1-\beta y_k)^{n-k},\cr
{\cal H}_{\rm id}(x,y)&=\prod_{k=1}^{n-1}(1+\beta x_k)^{n-k}.}
$$

\bigbreak
{\bf 2.4. Pieri's rule for Grothendieck polynomials.}
\medbreak

Let us denote by $(i,j)$ the transposition that interchanges $i$ and
$j$, $i<j$.

\medbreak
{\bf Proposition 5} ([FL]). {\it
$$G_{s_k}(x)\cdot G_w(x)\equiv\sum_v\beta^{l(v)-l(w)-1}G_v(x)~\mod I_n,
\eqno (4)
$$
where sum runs over $v\in S_n$, such that

$i)~~v=w\cdot (i_1,j_1)\ldots (i_{m+1},j_{m+1})$,

$ii)~~i_l\le k<j_l$ for all $l=1,2,\ldots ,m+1$,

$iii)~~ l(v)=l(w)+m+1$.}

More generally, consider the ring
$$P_{n,\beta }[y]:={\bf Z}[\beta ][x_1,\ldots ,x_n,y_1,\ldots ,y_n]
$$
and the ideal $J_n:=J_{n,\beta}$ in $R_{n,\beta}[y]$ generated by the following
polynomials 
$$e_i(x)+(-1)^{i-1}e_i(y),~~ i=1,\ldots ,n.
$$

\medbreak
{\bf Proposition 6} (Pieri's rule for double Grothendieck polynomials). {\it
$$G_{s_k}(x,y_w)\cdot G_w(x,y)\equiv G_{\rm id}(x,y_w)\left [\ds\sum_v
\beta^{l(v)-l(w)-1}G_v(x,y)\right ]~\mod J_n,\eqno (5)
$$
where $y_w=(y_{w(1)},\ldots ,y_{w(n)})$, and the sum runs over the same set 
of permutations $v\in S_n$ as in Proposition 5.}

\bigbreak
{\bf 2.5. Canonical involution $\omega$.}
\medbreak

There exists an involution $\omega$ of the ring $P_{n,\beta}[y]$ given by 
$\omega (x)={\buildrel\leftarrow\over x}$, $\omega (y)={\buildrel\leftarrow
\over y}$, where
for any sequence $z=(z_1,\ldots ,z_m)$ we define ${\buildrel\leftarrow\over
z}$ to be equal to $(z_m,z_{m-1},\ldots ,z_2,z_1)$. It is clear that 
involution $w$ preserves the ideals $I_n$ and $J_n$.

\medbreak
{\bf Proposition 7.} {\it Let $v$ be a permutation, then
$$\omega (G_v(x,y)){\cal H}_{\rm id}(x,y)\equiv (-1)^{l(v)}{\cal H}_{w_0vw_0}
(x,y)\omega (G_{\rm id}(x,y))~\mod J_{n,\beta}.\eqno (6)
$$
}

\medbreak
{\bf Remark.} Let $G_w^{LS}(A,B)$ be the double Grothendieck polynomials 
introduced by A.~Lascoux and M.P.-Sch\"utzenberger in [L], (2.1). Then
$$G_w^{LS}(1-\beta y,1+\beta x)=(-1)^{l(w)}{G_{w^{-1}}(x,y)\over G_{\rm id}
(x,y)}.
$$
\vfil\eject

\bigbreak
{\bf 2.6. Double Schubert polynomials.}
\medbreak

Let us remind the definition of double Schubert polynomials due to A.~Lascoux
and M.P.-Sch\"utzenberger.
\medbreak

{\bf Definition 2} (A.~Lascoux-- M.P.-Sch\"utzenberger). {\it For each 
permutation $w\in S_n$, the double Schubert polynomial $\s_w(x,y)$ is
defined to be
$$\s_w(x,y)=\partial_{w^{-1}w_0}^{(x)}\s_{w_0}(x,y),
$$
where divided difference operator $\partial_{w^{-1}w_0}^{(x)}$ acts on the $x$
variables, and
$$\s_{w_0}(x,y)=G_{w_0}(x,y)=\prod_{i+j\le n}(x_i+y_j).
$$
}

It follows from Definition 1 and 2 that
$$G_w(x,y)|_{\beta =0}={\cal H}_w(x,y)|_{\beta =0}=\s_w(x,y).
$$

Double Schubert polynomials appear in algebra and geometry as cohomology 
classes related to degeneracy loci of flagged vector bundles. If $h:E\to 
F$ is a map of rank $n$ vector bundles on a smooth variety $X$,
$$E_1\subset E_2\subset\cdots\subset E_n=E, \ F:=F_n\to 
F_{n-1}\to \cdots \to F_1
$$
are flags of subbundles and quotient bundles, then there is a degeneracy 
locus $\Omega_w(h)$ for each permutation $w$ in the symmetric group 
$S_n$, described by the conditions
$$\Omega_w(h)=\{ x\in X~|~{\rm rank}(E_p(x)\to F_q(x))\le\#\{ i\le q, 
w_i\le p\} , \forall p,q\}.
$$
For generic $h$,\  $\Omega_w(h)$ is irreducible, codim$\Omega_w(h)=l(w)$, 
and the class $[\Omega_w(h)]$ of this locus in the Chow ring of $X$ is 
equal to the double Schubert polynomial $\s_{w_0w}(x,-y)$, where
$$\eno{ 
x_i &= c_1(\ker (F_i\to F_{i-1})),\cr
y_i &= c_1(E_i/E_{i-1}), \ 1\le i\le n.}
$$
It is well-known [F1] that the Chow ring of flag variety $Fl_n$ admits the 
following description
$$CH^*(Fl_n)\cong{\bf Z}[x_1,\ldots ,x_n,y_1,\ldots ,y_n]/J,
$$
where $J$ is the ideal generated by
$$e_i(x_1,\ldots ,x_n)-e_i(y_1,\ldots ,y_n), \ 1\le i\le n,
$$
and $e_i(x)$ is the $i$-th elementary symmetric function in the variables 
$x_1,\ldots ,x_n$.

Result below describes the structure of quotiont ring ${\bf Z}[x,y]/J$.
\bigbreak

{\bf Proposition } ([LS], [KV]). {\it The ring ${\bf Z}[x_1,\ldots ,x_n,y_1,
\ldots ,y_n]/J_n$ is 
a free module of dimension $n!$ over the ring $R$, with basis either 
$\s_w(x)$, or $\s_w(x,y),\ w\in S_n$, where}
$$R:={{\bf Z}[x_1,\ldots ,x_n]\otimes{\bf Sym}[y_1,\ldots ,y_n]\over J}.
$$
\vfil\eject

\bigbreak
{\bf 2.7. Chern classes.}
\medbreak

Let us denote by $1+P_n^+$ the multiplicative monoid of polynomials in 
two sets of 
variables $x=(x_1,\ldots ,x_n)$ and $y=(y_1,\ldots ,y_n)$ with rational 
coefficients and constant term 1. There exists a homomorphism (the
Chern homomorphism)
$$c:{\bf Z}[x^{\pm 1},y^{\pm 1}]\to 1+P_n^+
$$
such that $c(1)=1$, $c(A+B)=c(A)c(B)$, $A,B\in{\bf Z}[x^{\pm},y^{\pm}]$.
On the basis $x^Iy^J$, $I\in{\bf Z}^n$, $J\in{\bf Z}^n$, it takes the 
values
$$c(x^Iy^J)=1+\sum_{k=1}^ni_kx_k-\sum_{k=1}^nj_ky_k.
$$

\medbreak
{\bf Proposition 8} (cf. [LS], \S 5).\ 
{\it Let us assume that $\beta =-1$. Then
$$\eno{
&\bullet c(y^{-\delta}G_v(1-x,y-1))\equiv 1-\sum_{v\le w}a_{vw}
\s_w(x,y)~\mod J_n,& (7)\cr
&\bullet c(x^{-\delta}{\cal H}_v(1-x,y-1))\equiv 1-\sum_{v\le w}b_{vw}
\s_w(x,y)~\mod J_n,}
$$
where $a_{vv}=b_{vv}=(-1)^{l(v)}(l(v)-1)!$.}

\vskip 0.4cm
{\bf \S 3. Quantization map.}
\bigbreak

In this section we describe a remarkable construction of quantization map,
which is due to S.~Fomin, S.~Gelfand and A.~Postnikov [FGP], see also [KM],
where some additional properties of quantization map are given. In [KM] the
quantization map was constructed using the Interpolation formula, [M], (6.8),
by replacing the classical double Schubert polynomials by their quantum 
analogs. Let us remind a 
construction from [FGP] in the form most convenient for our purposes. 
The starting point is a remarkable family of commuting operators
$$X_j:=X_j^{(n)}=x_j-\sum_{i<j}q_{ij}\partial_{(ij)}+\sum_{j<k}q_{jk}
\partial_{(jk)},~~1\le j\le n,
$$
where (for $i<j$) $\partial_{(ij)}=\partial_{t_{ij}}=\partial_i\partial_{i+1}
\ldots\partial_{j-1}\ldots\partial_{i+1}\partial_i$ is the divided difference
operator corresponding to the transposition $t_{ij}$, and $q_{ij}=q_iq_{i+1}
\ldots q_{j-1}$. It is clear that $X_i:P_n\to P_n$, but the following 
property is less obvious.
\bigbreak

{\bf Proposition} ([FGP], Theorem 3.1). {\it The operators $X_j$ commute 
pairwise.}
\medbreak

Another useful property of operators $X_i$ is the following.
\bigbreak

{\bf Lemma} ([FGP], Lemma 3.3). {\it For any polynomial $f\in P_n$, there 
exists a unique operator $F\in{\bf Z}[q_1,\ldots ,q_{n-1}][X_1,\ldots ,X_n]$
such that
$$F(X_1,\ldots ,X_n)(1)=f(x_1,\ldots ,x_n).
$$
}

The map $f\mapsto F:P_n\to P_n$ is called the quantization map, [FGP], 
Section~3. Proposition below describes one of the main property of 
quantization map.
\bigbreak

{\bf Proposition} ([FGP], (5.1)). {\it Let $\wt{\s}_w(x)$ be the quantum 
Schubert polynomial corresponding to a permutation $w\in S_n$. Then
$$\wt{\s}_w(X_1,\ldots ,X_n)(1)=\s_w(x_1,x_2,\ldots ,x_n).\eqno (9)
$$
}

The relation (9) is very useful and allows to reduce the proofs of 
statements about quantum Schubert polynomials to those about classical ones,
see e.g. [FGP], proof of Theorem~7.8. However, it seems rather difficult
to find an explicit (e.g. of Lascoux-Sch\"utzenberger's type) formula for
a quantum Schubert polynomial $\wt{\s}_w(x)$ using only the relation (9).
So we came to the problem how to describe effectively the quantization map.
To solve this problem, in the paper [FGP], Section~4, the basis of standard 
elementary monomials are used. However, this basis does not give an 
orthogonal basis with respect to the canonical pairing on the ring of
polynomials ${\bf Z}[x_1,\ldots ,x_n]$, see e.g. [LS]; [M], (5.2).
In the paper [KM] we gave another approach to the quantization problem which
is based on the theory of quantum double Schubert polynomials and quantum
Cauchy's identity. Our main observations are:

$i)$ with respect to the $y$ variables, the quantum double Schubert polynomials
have a structure which is very similar to that in the classical case;

$ii)$ the quantum Cauchy identity, in essential, is equivalent to the 
orthogonality of quantum Schubert polynomials with respect to quantum pairing,
see [KM], (2.5).

Quantum double Schubert polynomials ${\wt{\s}}_w(x,y)$,  
were defined in [KM] by using the Lascoux-Sch\"utzenberger type 
construction. Namely, by Definition 6, [KM],
$$\wt{\s}_w(x,y)=\partial_{w^{-1}w_0}^{(y)}\wt{\s}_{w_0}(x,y),
$$
where 
$$\wt{\s}_{w_0}(x,y)=\Delta_1(y_{n-1}|x_1)\Delta_2(y_{n-1}|x_1,x_2)
\ldots\Delta_{n-1}(y_1|x_1,\ldots ,x_{n-1})
$$ 
and $\Delta_k(t|x_1,\ldots ,x_k)$ is the Givental-Kim determinant [GK]:
$$\eno{
\Delta_k(t|x)&:=\det\pmatrix{x_1+t & q_1&0 &\ldots &\ldots 
&\ldots &0\cr
-1 & x_2+t & q_2 &0 &\ldots &\ldots & 0\cr
0 & -1 &x_3+t & q_3 & 0 &\ldots & 0 \cr
\vdots &\ddots &\ddots & \ddots &\ddots &\ddots &\vdots \cr
0&\ldots & 0 &-1&x_{k-2}+t &q_{k-2} & 0 \cr
0 &\ldots &\ldots & 0 &-1 & x_{k-1}+t & q_{k-1}\cr
0 &  \ldots &\ldots &\ldots & 0 & -1 & x_k+t}& (2)\cr
&=\sum_{i=1}^kt^{k-i}\wt e_i(x_1,\ldots ,x_k~|~q_1,\ldots ,q_{k-1}).}
$$ 
The polynomials $\wt e_i(x):=\wt e_i(x_1,\ldots ,x_k~|~q_1,\ldots ,q_{k-1}$
are called the quantum elementary symmetric polynomials degree $i$ in the
variables $x_1,\ldots ,x_k$. 

The main goal of this paper is to apply the results from [FGP] in the
context of [KM]. In sequel, we will use the notations from [KM]; e.g. we
denote by $\wt f$ the quantization of polynomial $f$. 

The following result will be used in \S 4.
\bigbreak

{\bf Theorem 1.} {\it Let $w_0$ be the longest element in $S_n$, then
$$\wt{\s}_{w_0}(X_1,\ldots ,X_n,y_1,\ldots , y_n)(1)=\s_{w_0}(x,y).
$$
}
{\it Proof.} We use two facts from [FGP], Corollary~4.6 and 
Proposition~4.10:

$\bullet$ $\wt e_i(X_1,\ldots ,X_n)(1)=e_i(x_1,\ldots ,x_n),$

$\bullet$ $\wt{fg}=f\wt g$, if $f\in\Lambda_n$.

Thus, we have \ $\Delta_{n-1}(y_1~|~X_1,\ldots ,X_{n-1})(1)=$
$$=\sum_ky_1^{n-1-k}\wt e_k(X_1,\ldots ,X_{n-1})(1)=
\sum_ky_1^{n-1-k}e_k(x_1,\ldots x_{n-1})=\prod_{k=1}^{n-1}(y_1+x_k).
$$
Hence,
$$\eno{
&\wt{\s}_{w_0}(X_1,\ldots ,X_{n-1},y_1,\ldots ,y_{n-1})(1)=
\Delta_1(y_{n-1}~|~X_1^{(n)})\cdots \Delta_{n-1}(y_1~|~X_1^{(n)},
\ldots ,X_{n-1}^{(n)})(1)\cr
&=\Delta_1(y_{n-1}~|~X_1^{(n-1)})\cdots\Delta_{n-2}(y_2~|~X_1^{(n-1)},
\ldots ,X_{n-2}^{(n-1)})(1)
\prod_{k=1}^{n-1}(y_1+x_k)=\prod_{i=1}^{n-1}\prod_{k=1}^{n-i}(y_i+x_k).}
$$
\qed

\vskip 0.4cm

{\bf \S 4. Quantum double and quantum double dual Grothendieck polynomials.}
\bigbreak

Let $x=(x_1,\ldots ,x_n)$, $y=(y_1,\ldots ,y_n)$ be two sets of variables,
and (cf. [KM])
$$\wt{\s}_{w_0}(x,y):=\wt{\s}_{w_0}^{(q)}(x,y)=\prod_{i=1}^{n-1}
\Delta_i(y_{n-i}~|~x_1,\ldots ,x_i),
$$
where $\Delta_k(t|x_1,\ldots ,x_k):=\sum_{j=0}^kt^{k-j}e_j(x_1,\ldots ,
x_k~|~q_1,\ldots ,q_{k-1})$ is the Givental-Kim determinant (2), 
which is generating function for the quantum 
elementary symmetric functions in the variables $x_1,\ldots ,x_k$. 

\bigbreak
{\bf Definition 3.} {\it For each permutation $w\in S_n$, the quantum double 
dual Grothendieck polynomial $\wt{\cal H}_w(x,y)$ is defined to be
$$\wt{\cal H}_w(x,y)=(\pi^-_{ww_0})^{(y)}\wt{\s}_{w_0}(x,y).
$$
}
Here symbol $(\pi^-_{ww_0})^{(y)}$ means that isobaric operator 
$\pi^-_{ww_0}$ acts on the $y$ variables.

\bigbreak
{\bf Corollary 1} (of Theorem 1). \ {\it for each permutation $w\in S_n$,}
$$\eno{
&\wt{\s}_w(X,y)(1)=\s_w(x,y),& (10)\cr
&\wt{\cal H}_w(X,y)(1)={\cal H}_w(x,y).}
$$
\bigbreak
{\bf Remarks.} 1) In the ``classical limit'' $q_1=\cdots =q_{n-1}=0$,
$$\wt{\cal H}_w(x,y)|_{q=0}=(\pi^-_{ww_0})^{(y)}G_{w_0}(x,y)=
G_{w^{-1}}^{(-\beta )}(y,x)={\cal H}_w^{(\beta )}(x,y),
$$
i.e. ${\cal H}_w(x,y)|_{q=0}={\cal H}_w^{(\beta )}(x,y)$.

2) (Stability) Let $m>n$ and let $i:S_n\hookrightarrow S_m$ be the 
embedding. Then
$$\wt{\cal H}_w(x,y)=\wt{\cal H}_{i(w)}(x,y).
$$

\vskip 0.3cm
Let us define the quantum dual Grothendieck polynomial $\wt{\cal H}_w(x)$
as the specialization
$$\wt{\cal H}_w(x):=\wt{\cal H}_w(x,y)|_{y=0}.
$$

\bigbreak
{\bf Definition 4.} {\it For each permutation $w\in S_n$, the quantum 
double Grothendieck polynomial $\wt G_w(x,y)$ is defined to be}
$$\wt G_w(x,y)=\sum_{w\le v}(-\beta )^{l(v)-l(w)}\wt{\cal H}_w(x,y).
$$

Let us define the quantum Grothendieck polynomials $\wt G_w(x)$ as
the specialization
$$\wt G_w(x):=\wt G_w(x,y)|_{y=0}.
$$

\bigbreak
{\bf Example.} Quantum double and quantum double dual Grothendieck 
polynomials for $S_3$.
$$\eno{
&\wt G_{121}(x,y)=(x_1+y_1)(x_1+y_2)(x_2+y_1)+q_1(x_1+y_2)\cr
&\wt G_{12}(x,y)=[(x_1+y_1)(x_2+y_1)+q_1](1-\beta y_2),\cr
&\wt G_{21}(x,y)=[(x_1+y_1)(x_1+y_2)-q_1](1-\beta y_1),\cr
&\wt G_1(x,y)=(x_1+y_1)[(1-\beta y_1)(1-\beta y_2)+q_1\beta ^2],\cr
&\wt G_2(x,y)=(x_1+x_2+y_1+y_2-\beta y_1y_2+\beta x_1x_2+\beta q_1)
(1-\beta y_1),\cr
&\wt G_{\rm id}(x,y)=(1-\beta y_1)[(1-\beta y_1)(1-\beta y_2)+q_1\beta^2];\cr
&\wt{\cal H}_{121}(x,y)=\wt G_{121}(x,y),\cr
&\wt{\cal H}_{12}(x,y)=[(x_1+y_1)(x_2+y_1)+q_1](1+\beta x_1),\cr
&\wt{\cal H}_{21}(x,y)=[(x_1+y_1)(X_1+y_2)-q_1](1+\beta x_2)+
q_1\beta (x_1+x_2+y_1+y_2),\cr
&\wt{\cal H}_1=(x_1+y_1)[(1+\beta x_1)(1+\beta x_2)+q_1\beta^2],\cr
&\wt{\cal H}_2=(x_1+x_2+y_1+y_2+\beta x_1x_2-\beta y_1y_2+\beta q_1)
(1+\beta x_1),\cr
&\wt{\cal H}_{\rm id}=(1+\beta x_1)[(1+\beta x_1)(1+\beta x_2)+q_1\beta^2].}
$$

{\bf Remark.} One can show that
$$\eno{
&(\pi_{w_0}^-)^{(y)}\wt G_{w_0}(x,y):=\wt G_{\rm id}(x,y)=
\beta^{n(n-1)\over 2}
\wt G_{w_0}\left( x,(\beta^{-1},\beta^{-1},\ldots ,\beta^{-1})\right),\cr
&\wt{\cal H}_{\rm id}^{(\beta )}(x,y)=\wt G_{\rm id}^{(-\beta )}(y,x).}
$$

{\bf Theorem 2.} (Quantum Cauchy identity for Grothendieck polynomials).
$$\sum_{w\in S_n}\wt{\cal H}_w(x,\ominus z)G_{ww_0}(y,z)=\prod_{k=1}^{n-1}
\Delta_k^{(\beta )}(y_{n-k}~|~x_1,\ldots ,x_k),\eqno (11)
$$
{\it where} $\Delta_k^{(\beta )}(t~|~x_1,\ldots ,x_k)=(1+\beta t)^k
\Delta_k\left(\ds{t\over 1+\beta t}~\mid ~x_1,\ldots ,x_k\right)$.

Proof follows from the classical Cauchy identity for Grothendieck 
polynomials (2), Theorem~1 and (10).

Let us denote the product in the RHS(11) by $\wt{\bf G}_{w_0}(x,y)$.

\medbreak
{\bf Corollary 2} (of Theorem 2).

$\bullet$ \ $ \wt{\cal H}_w(x)=\psi_{ww_0}^{(y)}\wt{\bf G}_{w_0}(x,y)|_{y=0}$,

$\bullet$ \ $ \wt G_w(x)=\pi_{ww_0}^{(y)}\wt{\bf G}_{w_0}(x,y)|_{y=0},$

$\bullet$ \ $ \wt G_w(x)=\sum_{v\in S_n}\eta (\pi_{ww_0}G_{vw_0}(x))
\wt{\cal H}_v(x).$

It seems interesting to study the properties of polynomials 
$\wt{\bf G}_w(x,y)=\pi_{ww_0}^{(y)}\wt{\bf G}_{w_0}(x,y)$ and 
$\wt{\bf H}_w(x,y)=\psi_{ww_0}^{(y)}\wt{\bf G}_{w_0}(x,y)$.

\vfil\eject
\vskip 1cm
{\bf References.}
\vskip 0.5cm

\item {[C]} Ciocan--Fontanine I., {\it Quantum cohomology of flag varieties},
Intern. Math.Research Notes, 1995, n.6, p.263-277;

\item {[F1]} Fulton W., {\it Flags, Schubert polynomials, degeneracy loci, and
determinantal formulas,} Duke Math. J., 1991, v.65, p.381-420;

\item {[F2]} Fulton W., {\it Young tableaux with applications to representation
theory and geometry,} Preprint, 1995;

\item {[FL]} Fulton W. and A. Lascoux A., {\it A Pieri formula in the
Grothendieck ring of a flag bundle,} Duke Math. J., 1994, v.76, p.711-729;

\item {[FGP]} Fomin S., Gelfand S. and Postnikov A., {\it Quantum Schubert
polynomials,} Preprint, 1996;

\item {[FK]} Fomin S. and Kirillov A.N., {\it Grothendieck polynomials and
the Yang-Baxter equation}, Proceedings of the 6th International Confrence
on Formal Power Series and Algebraic Combinatorics, DIMACS, 1994, p.183-190;

\item {[G]} Grothendieck A., {\it Theoreme de Riemann--Roch,} in SGA 6,
Springer L.N., (1971), v.225;

\item {[GK]} Givental A. and Kim B., {\it Quantum cohomology of flag manifolds 
and Toda lattices}, Comm.Math.Phys., 1995, v.168, p.609-641;    

\item {[K]} Kim B., {\it Quot schemes for flags and Gromov invariants for flag 
varieties}, Prepr., 1995, alg-geom/9512003;

\item {[KM]} Kirillov A.N. and Maeno T., {\it Quantum Schubert polynomials and
Vafa-Intriligator formula,} Preprint UTMS 96-41, 1996, 50p.;

\item {[L]} Lascoux A., {\it Anneau de Grothendieck de la variete de 
drapeaux,} in ``The Grothendieck Festtchift'', vol III, Birkhauser, 
1990, p.1-34;

\item {[LS]} Lascoux A. and Sch\"utzenberger M.-P., {\it Symmetry and flag
manifolds,} Lect. Note in Math., 1983, v.996, p.118-144;

\item {[M]} Macdonald I.G., {\it Notes on Schubert polynomials}, Publ. LCIM, 
1991, Univ. de Quebec a Montreal.

\end